\documentclass{amsart}

\usepackage{amsmath,amssymb}
\usepackage{graphicx}
\usepackage{enumerate}

\makeatletter
\@addtoreset{equation}{section}
\makeatother

\numberwithin{equation}{section}

\def\res{\mathrm {res}}

\def\Sym{\mathrm {Sym}}

\def\CC{{\mathbb C}}

\newtheorem{definition}{Definition}[section]

\newtheorem{theorem}{Theorem}[section]
\newtheorem{proposition}{Proposition}[section]

\newtheorem{lemma}{Lemma}[section]

\def\book#1{\rm{#1}, }
\def\paper#1{\textit{#1}, }
\def\jour#1{\rm{#1}, }
\def\yr#1{({\rm{#1}) }}
\def\vol#1{\textbf{#1}}
\def\pages#1{\rm{#1}}

\def\publaddr#1{\rm{#1}, }
\def\publ#1{\rm{#1}, }
\def\by#1{{\rm{#1}, }}

\def\eds{\rm{eds.}}

\begin{document}

\title{$UVW$ relations over
a subvariety of a hyperelliptic
Jacobian}

\author{Shigeki MATSUTANI}
\address{8-21-1 Higashi-Linkan, Sagamihara, 228-0811, JAPAN}

\email{rxb01142@nifty.com}

\maketitle

\begin{abstract}
This article extends relations of Mumford's $UVW$-expressions
to those in subvariety in a hyperelliptic Jacobian
using Baker's method.
\end{abstract}



\section{ Introduction}

For a hyperelliptic curve $C_g$ whose affine part is
given by $y^2 = \prod_{i=1}^{2g+1}(x-b_i)$,
where $b_i$'s are complex numbers,
its Jacobian $\mathcal J_g$ is given as a complex
torus $\CC^g/\Lambda$ by the Abel map $\omega$ \cite{Mu}.
The Abelian theorem enables us to have a natural
morphism from the symmetrical product $\Sym^g(C_g)$
to the Jacobian
$\mathcal J_g\approx\omega[\Sym^{g}(C_g)]/\Lambda$.

Mumford and his coworkers used $UVW$
expression based upon Jacobi's considerations \cite{Mu},
which  represents the hyperelliptic functions over the
Jacobian

Let $D(z)$ be a certain derivative of the Jacobian,
the Bolza polynomial, $F(z)\equiv U(z):=(z-x_1)\cdots(z-x_g)$
for $(x_i, y_i)_{i=1, \cdots, g} \in \Sym^g (C_g)$.
Further let
$V(z) := D(z) (x_1 + \cdots + x_2)$ and
$W(z) := (f(z) - V(z)^2)/F(z)$.
Mumford and his coworkers showed [\cite{Mu} Theorem 3.1],
\begin{theorem}\label{th:1}
\begin{enumerate}
\item
$$
	D(z_1) F(z_2) = \frac{F(z_2) V(z_1) - V(z_2) F(z_1)}
	                  {z_1- z_2}.
$$
\item
$$
	D(z_1) V(z_2) = \frac{1}{2}\left(\frac{F(z_2) W(z_1) - W(z_2) F(z_1)}
	                  {z_1- z_2} -F(z_1) F(z_2)\right),
$$

\item
$$
	D(z_1) W(z_2) = \frac{W(z_2) V(z_1) - V(z_2) W(z_1)}
	                  {z_1- z_2} +F(z_1) V(z_2).
$$

\item $D(z_1)D(z_2)=D(z_2)D(z_1)$.

\end{enumerate}
\end{theorem}
There exists an interesting Poisson structure in these
relations, which are studied by several authors
\cite{AHP, PV, Mu}.

On the other hand, as zeros of an appropriate shifted Riemann theta
function over $\mathcal J_g$, the theta divisor
is defined as
$$
      \Theta := \omega[\Sym^{g-1}(C_g)]/\Lambda
$$
which is a subvariety of $\mathcal J_g$.
Similarly, it is natural that we introduce a subvariety
$$
      \Theta_k := \omega[\Sym^{k}(C_g)]/\Lambda
$$
and a sequence,
$$
    \Theta_0\subset\Theta_1 \subset \Theta_2 \subset \cdots \subset
    \Theta_{g-1} \subset \Theta_g \equiv \mathcal J_g
$$
Vanhaecke studied the structure of these subvarieties
as stratifications of the Jacobian
 $\mathcal J_g$ using the strategies
developed the studies in the  integrable system \cite{V, V2}.
He showed that it is connected with  stratifications of the
Sato Grassmannian \cite{V}.
Further in \cite{V2},
he studied Lie-Poisson structure in the
Jacobian and showed that  invariant
manifolds associated with Poisson brackets
can be identified with
these strata; it implies that these strata
are characterized by the Lie-Poisson structure.
Further he and, Abenda and Fedorov \cite{AF}
considered their relations
 to  finite-dimensional integrable systems, {\it i.e.},
 Henon-Heiles system and Neumann systems.
Independently the author considered a relation of
symmetric functions over $\Theta_k$ as
an extension of the study of Weierstrass on
al-functions \cite{Ma}.

The elementary symmetric functions over $\Theta_k$
appear in \cite{AF, Ma} and play the important roles
to reveal structure of $\Theta_k$.
In the case of the Jacobian, the relations of
the elementary symmetric functions over the Jacobian
is represented by Theorem \ref{th:1},
which
is directly related to Neumann system and other many
studies on structures, like a Lie-Poisson structure,
 of hyperelliptic curves \cite{AHP, PV, Mu}.
Though the structure of these subvarieties
was studied using Theorem \ref{th:1} \cite{V2},
its variant over $\Theta_k$ was not studied.

Thus the purpose of this article is an extension
of the relations in Theorem \ref{th:1} to
similar variants  of elementary symmetric functions over
$\Theta_k$ as in our main theorem \ref{th:3}.
We believe
that Theorem \ref{th:3} has an effects on these studies.

Of course, in this stage our theorem \ref{th:3} is not
connected with such a finite integrable system directly
 though
Abenda and Fedorov \cite{AF} studied similar subjects,
and E. Previato suggested the author that
there might be
a connection between Theorem \ref{th:3}
and a finite integrable system.
We expect that our results shed some light
on these studies.

Furthermore our strategy in this article is based upon
Baker's method in \cite{Ba}, which is a direct application
of the reciprocity laws for differentials over
a curve to a relation over there. Thus we believe
that we should re-evaluate Baker's
method using modern expressions of the reciprocity
\cite{BP} in future.
If we could, it is expected that
we would have a modern language to expresses
the subvarieties in Jacobians.

\bigskip

A final step in this work was
done at Concordia University and
thus the author thanks Professor J. McKay for his kindness
and hospitality.
The author is sincerely  grateful to Professor
E. Previato for her kind comments
on this work and fruitful discussions.
He also thanks Professor Y. \^Onishi for his continuous supports
on hiss study.

\bigskip
\section{Hyperelliptic curve}

\textit{ Hyperelliptic Curve:}
\rm{
This article deals with a hyperelliptic curve $C_g$  of genus $g$
$(g>0)$ given by the affine equation,
\begin{gather*} \begin{split}
   y^2 &= f(x) \\  &= \lambda_{2g+1} x^{2g+1} +
\lambda_{2g} x^{2g}+\cdots  +\lambda_2 x^2
+\lambda_1 x+\lambda_0  \\
     &=b_0 (x-b_1)(x-b_2)\cdots (x-b_{2g+1}),\\
\end{split}  \label{2-1}
\end{gather*}
where
 $\lambda_j$'s
and $b_j$'s  are  complex numbers.}

For a point $(x_i, y_i)\in C_g$,
the unnormalized differentials of the first kind are
defined by,
\begin{gather*}   d u^{(i)}_1 := \frac{ d x_i}{2y_i}, \quad
      d u^{(i)}_2 :=  \frac{x_i d x_i}{2y_i}, \quad \cdots, \quad
     d u^{(i)}_g :=\frac{x_i^{g-1} d x_i}{2 y_i}.
      \label{2-3}
\end{gather*}

For positive integer $k(\leq g$.
the Abel map from $k$-th symmetric product
of the curve $C_g$ to $\CC^g$ is defined by,
\begin{gather*}
u :=(u_1,\cdots,u_g)
:\mathrm{Sym}^k( C_g) \longrightarrow \mathbb C^g,
\quad
  u_k:= \sum_{i=1}^k
       \int_\infty^{(x_i,y_i)} d u^{(i)}_k.
      \label{2-6}
\end{gather*}
The Jacobian is defined by
$$
	\mathcal J_g := \mathbb C^g / <\mbox{lattice}>.
$$

Let its image quotient by the lattice denoted by
$\Theta_k$.
\begin{gather}
 \{0\} \subset C_g \equiv \Theta_1 \subset \cdots
  \subset \Theta_{g-1} \subset \Theta_{g} \equiv \mathcal J_g.
\end{gather}

\bigskip

Let us fix $k\le g$ and
 $((x_1,y_1), (x_2,y_2), \cdots, (x_k,y_k))\in\mathrm{Sym}^k( C_g)$.
We will introduce a variant  in $\Theta_k$ of $UVW$-expression,
which also appears in \cite{AF, Ma, Mu}.
\begin{definition}
We define
\begin{enumerate}[{(}1{)}]

\item
$
F^{(k)}(z) :=  (z-x_1)(z-x_2)\cdots (z-x_k)
$
and for brevity we denote it by $F(z)$ if there is
no confusion.

\item Let $k$ be an integer such that $k \le m \le g$
and natural number $n := m - k + 1$.
$\displaystyle{
D(z):=\frac{1}{2z^n} \sum_{i=n}^m z^{i}\partial_{u_i}}
$

\item
$
V(z) := D(z) (x_1+x_2+\cdots +x_k).
$

\item
$
W(z) := (f(z)-[V(z)z^{n-1}]^2 )/F(z)z^{2n-2}.
$
\end{enumerate}
\end{definition}

Simple consideration gives the following proposition:
\begin{proposition} \label{prop:2-1}
\begin{enumerate}
\item
$$
   D(z) = \sum_{i=1}^k \frac{y_i F(z)}{F'(x_i)(z-x_i)x_i^n}
                       \frac{\partial}{\partial x_i},
$$

\item
 $$
	V(z) = \sum_{i=1}^k \frac{y_i F(z)}{F'(x_i)(z-x_i)x_i^n}.
$$

\item $y_i=V(x_i)x_i^n$.
\end{enumerate}
\end{proposition}

\begin{proof}
Let us introduce quantities and matrices:
\begin{gather*}
\begin{split}
\pi_i(x) &:= \frac{F(x)}{x-x_i}\\
        &=\chi_{i,k-1}x^{k-1} +\chi_{i,k-2} x^{k-2}
            +\cdots+\chi_{i,1}x+\chi_{i,0},\\
\end{split}
\end{gather*}
\begin{gather*}
{\mathcal W} := \begin{pmatrix}
     \chi_{1,0} & \chi_{1,1} & \cdots & \chi_{1,k-1}  \\
      \chi_{2,0} & \chi_{2,1} & \cdots & \chi_{2,k-1}  \\
   \vdots & \vdots & \ddots & \vdots  \\
    \chi_{k,0} & \chi_{k,1} & \cdots & \chi_{k,k-1}
     \end{pmatrix},\quad
	\mathcal Y := \begin{pmatrix}
     y_1/x_1^{n-1} & \ & \ & \  \\
      \ & y_2/x_2^{n-1}& \ & \   \\
      \ & \ & \ddots   & \   \\
      \ & \ & \ & y_k /x_k^{n-1} \end{pmatrix},
\end{gather*}
\begin{gather*}
	\mathcal F' := \begin{pmatrix}
     F'(x_1) & \ & \ & \  \\
      \ & F'(x_2)& \ & \   \\
      \ & \ & \ddots   & \   \\
      \ & \ & \ & F'(x_{k})  \end{pmatrix},
\quad
{\mathcal V}:= \begin{pmatrix} 1 & 1 & \cdots & 1 \\
                   x_1 & x_2 & \cdots & x_k \\
                   x_1^2 & x_2^2 & \cdots & x_k^2 \\
                    \vdots& \vdots &       & \vdots \\
                   x_1^{k-1} & x_2^{k-1} & \cdots & x_k^{k-1}
                 \end{pmatrix},
\end{gather*}
where $F'(x):=d F(x)/d x$.  Then we have
\begin{gather*}
\begin{pmatrix}d u_n\\ \vdots\\ d u_m\end{pmatrix}
= \frac{1}{2}{\mathcal V}{\mathcal Y}^{-1}
\begin{pmatrix}d x_1\\ \vdots\\ d x_k\end{pmatrix}.
\end{gather*}
and ${\mathcal V}^{-1}= {\mathcal F'}^{-1}{\mathcal W}$.
By letting $\partial_{u_i}:=\partial/\partial{u_i}$ and
$\partial_{x_i}:=\partial/\partial{x_i}$, we obtain
\begin{gather*}
	\begin{pmatrix} \partial_{u_n}\\
                 \partial_{u_{n+1}}\\
                 \vdots\\
                 \partial_{u_m}
         \end{pmatrix}
   = {}^t(2\mathcal Y \mathcal F^{\prime -1}\cdot {\mathcal W})
	\begin{pmatrix} \partial_{x_1}\\
                 \partial_{x_2}\\
                 \vdots\\
                 \partial_{x_k}
         \end{pmatrix}.
\end{gather*}
Hence $D(z)$ is given by (1).
(2) and (3) are obvious from (1).
\end{proof}

\bigskip
\section{Relations of $UVW$ in $\Theta_k$}

Now we will give our main theorem as follows:

\begin{theorem}\label{th:3}
Assume $\lambda_0\neq0$.
\begin{enumerate}
\item
$$
	D(z_1) F(z_2) = \frac{F(z_2) V(z_1) - V(z_2) F(z_1)}
	                  {z_1- z_2}.
$$
\item
$$
	D(z_1) V(z_2) = \frac{1}{2}\left[\frac{F(z_2) W(z_1) - W(z_2) F(z_1)}
	                  {z_1- z_2} -F(z_1) F(z_2)H(\lambda, z_1, z_2)\right],
$$
where $H(\lambda, z_1, z_2):=H_0(\lambda, z_1, z_2)
+H_\infty(\lambda, z_1, z_2)$,
\begin{gather*}
\begin{split}
	H_0(\lambda, z_1, z_2)&:=-
  \frac{1}{(2n-3)!}
      \left[\frac{\partial^{2n-3}}{\partial x^{2n-3}}
         \frac{f(x)}{F(x)^2(z_1 - x)(z_2- x)}\right]_{x=0},\\
H_\infty(\lambda, z_1, z_2)&:=\left(
    \frac{1}{(4g-4m)!}
     \left[\frac{\partial^{4g-4m}}{\partial t^{4g-4m}}
         \frac{t^{4g- 4k-2}f(1/t^2)}
        {F(1/t^2)^2(1 - z_1 t^2)(1 - z_2 t^2)}\right]_{t=0}\right).
\end{split}
\end{gather*}
\item
$$
	D(z_1) W(z_2) = \frac{W(z_2) V(z_1) - V(z_2) W(z_1)}
	                  {z_1- z_2} +F(z_1) V(z_2)H(\lambda, z_1, z_2).
$$

\item $D(z_1)D(z_2)=D(z_2)D(z_1)$.

\end{enumerate}
\end{theorem}

From here, we will give the proof of this theorem.
Its idea is simple. After translating the words in $u$'s into
$(x_i, y_i)$'s, we will perform a residual computations
around a boundary of a polygon representation
$C_g^o$ of $C_g$, and use some combinatorial trick in sums.

\begin{proof}
Using Proposition \ref{prop:2-1},
(1) is directly shown: $ D(z_1) F(z_2)$ is equal to
\begin{gather*}
\begin{split}
   & -\sum_{i=1}^k
   \frac{F(z_1)F(z_2) y_i }{F'(x_i)x_i^{n-1}(z_1-x_i)(z_2-x_i)}\\
   &=-\sum_{i=1}^k
   \left(\frac{F(z_1)F(z_2) y_i}{F'(x_i)x_i^{n-1}(z_1-x_i)}
  -
   \frac{F(z_1)F(z_2) y_i }{F'(x_i)x_i^{n-1}(z_2-x_i)}\right)
   \frac{1}{z_1-z_2},
 \end{split}
\end{gather*}
which  is the right hand side of (1) from the definition of $V$.

Provided (2), (3) is proved as follows:
Noting that $D(z_1)$ is a direction differential
operator and its action obeys the Leibniz rule.
Using the fact $D(z_1)f(z_2)=0$,
and the Leibniz rule, $D(z_1)W(z_2)$ is given by
\begin{gather*}
  \frac{-1}{F(z_2)^2z_2^{2n-2}}
[D(z_1)F(z_2)] (f(z_2)-V(z_2)^2 z_2^{2n-2})
-\frac{2}{F(z_2)}[V(z_2) D(z_1)V(z_2)].
\end{gather*}
Using (1) and (2), it becomes
\begin{gather*}
\begin{split}
& \frac{-1}{F(z_2)^2z_2^{2n-2}}(f(z_2)-V(z_2)^2 z_2^{2n-2})
\frac{V(z_2) F(z_1) - F(z_2) V(z_1)}
	                  {z_1- z_2}\\
&-	\frac{2}{F(z_2)}V(z_2)
\left(\frac{1}{2}\frac{F(z_2) W(z_1) - W(z_2) F(z_1)}
	                  {z_1- z_2} -\frac{1}{2}F(z_1) F(z_2)H(\lambda, z_1, z_2)
	                  \right),
\end{split}
\end{gather*}
which  is the right hand side of  (3) from the definition of $W$.

(4) is roughly proved due to the fact
$D(z_1) V(z_2) =D(z_2) V(z_1)$.
Thus we consider (4) after proving (2).

Let us consider the formula (2).
The strategy is essentially the same as \cite{Ba}.
First we translate the words of the Jacobian into those of
the curves in $\mathrm{Sym}^k(C_g)$;
we rewrite the differentials in the Jacobian
in terms of the differentials over curves
 in $\mathrm{Sym}^k(C_g)$ by (1). We count the
residue of an integration and use a combinatorial trick.
Then we will obtain (2):
\begin{gather*}
D(z_1)V(z_2) =
 \sum_{j=1, i=1}^k \frac{F(z_1)y_i}{F'(x_i)x_i^{n-1}(z_1-x_i)}
        \frac{\partial}{\partial x_i}
         \frac{F(z_2)y_j}{F'(x_j)x_j^{n-1} (z_2-x_j)}.
\end{gather*}
Let us decompose the summation into $j=i$ and $j\neq i$ parts.
Further we will note the derivative of $F(x)$;
\begin{gather*}
	\frac{\partial}{\partial x_i} \left(
   \left[\frac{\partial}{\partial x} F(x)\right]_{x=x_i}\right)
         =\frac{1}{2}
 \left[\frac{\partial^2}{\partial x^2} F(x)\right]_{x=x_i}
         .
\end{gather*}
Then we obtain the formulation that $D(z_1)V(z_2)$ is equal to
\begin{gather*}
\begin{split}
 &
 \sum_{ i=1}^k \frac{F(z_1)F(z_2)}{F'(x_i)x_i^{2n-2}(z_2-x_i)(z_1-x_i)}
 \frac{1}{2}\left(
\frac{f'(x_i)}{F'(x_i)}-
\frac{f(x_i)F''(x_i)}{F'(x_i)}
-2(n-1)\frac{f(x_i)}{x_i F'(x_i)}\right)\\
&
+
 \sum_{i\neq j}^k \frac{F(z_1)F(z_2)y_i y_j}
{F'(x_i)F'(x_j)x_i^{n-1}x_j^{n-1}(z_2-x_i)(z_1-x_j)}
\left(\frac{1}{x_j-x_i}-\frac{1}{z_2 - x_i}\right)\\
\end{split}
\end{gather*}
\begin{gather*}
\begin{split}
=
\frac{F(z_1)F(z_2)}{2} \sum_{ i=1}^k \frac{1}{F'(x_i)}
\left(\frac{\partial}{\partial x}
\frac{f(x)}{F'(x)^2x^{2n-2}(x-z_1)(x-z_2)}\right)_{x=x_i}\\
+\sum_{ i=1}^k \frac{F(z_1)F(z_2)}{F'(x_i)x_i^{2n-2}(z_2-x_i)(z_1-x_i)}
 \frac{1}{2}\left(
\frac{1}{z_2-x_i}+\frac{1}{z_1-x_i}\right)\\
  +
 \sum_{i\neq j}^k \frac{F(z_1)F(z_2)y_i y_j}
{F'(x_i)F'(x_j)x_i^{n-1}x_j^{n-1}(z_2-x_i)(z_1-x_j)}
\left(\frac{1}{x_j-x_i}-\frac{1}{z_2 - x_i}\right).\\
\end{split}
\end{gather*}
We will consider each term in the formulae; we refer the first
and second terms
$[DV]_1(z_1,z_2)$ and the third term $[DV]_2(z_1,z_2)$.
The proof of (2) finishes due to the following Lemma because
$[DV]_1(z_1,z_2)$ is given by
\begin{gather*}
\begin{split}
 &
\frac{1}{2}\frac{F(z_2) f(z_1)/z_2^{2n-2} F(z_1)
- f(z_2) F(z_1)/z_1^{2n-2} F(z_2)}
	                  {z_1- z_2} -F(z_1) F(z_2)H(\lambda, z_1, z_2)\\
 &+\frac{1}{2}\frac{1}{z_1-z_2}
\sum_{ i=1}^k \frac{F(z_1)F(z_2)}{F'(x_i)x_i^{2n-2}}
\left(\frac{1}{(z_2-x_i)^2}-\frac{1}{(z_1-x_i)^2}\right)\\
\end{split}
\end{gather*}
whereas $[DV]_2(z_1,z_2)$ is equal to
\begin{gather*}
 \frac{1}{2}\frac{1}{z_1-z_2}
\sum_{i\neq j}^k \frac{F(z_1)F(z_2)y_i y_j}
{F'(x_i)F'(x_j)x_i^{n-1}x_j^{n-1}}\left(\frac{1}{(z_1-x_j)(z_1-x_i)}
-\frac{1}{(z_2 - x_i)(z_2-x_j)}\right).
\end{gather*}
The first term in  $[DV]_1(z_1,z_2)$ is equal to  parts of $W$
in the first term in (2) and the second term in (2).
The remainder is given by the second term in  $[DV]_1(z_1,z_2)$
and $[DV]_2(z_1,z_2)$, which is
\begin{gather*}
 \frac{1}{2}\frac{1}{z_1-z_2}
\sum_{i, j}^k \frac{F(z_1)F(z_2)y_i y_j}
{F'(x_i)F'(x_j)x_i^{n-1}x_j^{n-1}}\left(\frac{1}{(z_1-x_j)(z_1-x_i)}
-\frac{1}{(z_2 - x_i)(z_2-x_j)}\right).
\end{gather*}
It is obvious that they forms the right hand side of (2).

Finally we will consider the commutativity of $D(z_a)$.
Since $D(z_1)D(z_2)$ is decomposed to the form
$$
     \sum_{i=1}^k [DD]_{i}(z_1, z_2)\frac{\partial}{\partial x_i}
     +\sum_{i,j=1}^k [DD]_{ij}(z_1, z_2)\frac{\partial^2}{
     \partial x_i \partial x_j},
$$
and obviously $[DD]_{ij}(z_1, z_2)=[DD]_{ij}(z_2, z_1)$
we must check the first term.
Let $[DV]_1(z_1,z_2)$ and $[DV]_2(z_1,z_2)$ be
denoted by
$$
	[DV]_1(z_1,z_2)=\sum_{i=1}^k [DV]_{1,i}(z_1,z_2),
\quad
	[DV]_2(z_1,z_2)=\sum_{i\neq j}^k [DV]_{2,ij}(z_1,z_2),
$$
following the definition in straightforward way.
Then we have
$$
 \sum_{i=1}^k [DD]_{i}(z_1, z_2)\frac{\partial}{\partial x_i}
 =\sum_{i=1}^k [DV]_{1,i}(z_1,z_2)
 \frac{\partial}{\partial x_i}
 +\sum_{i\neq j}^k [DV]_{2,ij}(z_1,z_2)
 \frac{\partial}{\partial x_j}.
 $$
Let us consider $D(z_1)D(z_2)-D(z_2)D(z_1)$.
Since it is clear that $[DV]_{1,i}(z_1,z_2)=[DV]_{1,i}(z_2,z_1)$,
$$
\sum_{i\neq j}^k [DV]_{2,ij}(z_1,z_2)
 \frac{\partial}{\partial x_j}-
 \sum_{i\neq j}^k [DV]_{2,ij}(z_2,z_1)
 \frac{\partial}{\partial x_j}
$$
remains. Using the appropriate symmetric quantities
$K_{2,ij}(z_1,z_2)=K_{2,ij}(z_2,z_1)=K_{2,ji}(z_1,z_2)$,
it becomes
\begin{gather*}
\begin{split}
\sum_{i\neq j}^k K_{2,ij}(z_1,z_2)
&\Bigr[
\frac{1}{z_1-x_i}
\frac{1}{z_2 -x_j}
\left(\frac{1}{x_j-x_i}-\frac{1}{z_2 - x_i}\right)
\\
&-
\frac{1}{z_2-x_i}
\frac{1}{z_1 -x_j}
\left(\frac{1}{x_j-x_i}-\frac{1}{z_1 - x_i}\right)\Bigr]
 \frac{\partial}{\partial x_j}.
\end{split}
\end{gather*}
Noting
$$
\frac{1}{z_1-x_i}
\frac{1}{z_2 -x_j}
\left(\frac{1}{x_j-x_i}-\frac{1}{z_2 - x_i}\right)
=\frac{1}{z_1-x_i}
\frac{1}{z_2 -x_i}\frac{1}{x_j-x_i},
$$
this vanishes. Hence (4) is proved.
\end{proof}

\begin{lemma}
Following relations hold:
\begin{enumerate}[{(}1{)}]
\item
\begin{gather*}
\begin{split}
\sum_{i=1}^k \frac{1}{F'(x_i)}
           &\left[\frac{\partial}{\partial x}\left(
       \frac{f(x)}{(x - z_1)(x-z_2)x^{2n-2} F'(x)} \right) \right]_{x = x_i}\\
      &= \frac{1}{(4g-4m)!}\left[
     \frac{\partial^{4g-4m}}{\partial t^{4g-4m}}
         \frac{t^{4g- 4k-2}f(1/t^2)}
        {F(1/t^2)^2(1 - z_1t^2)(1 - z_2 t^2)}\right]_{t=0}\\
      &  -\frac{1}{(2n-3)!}
      \left[\frac{\partial^{2n-3}}{\partial x^{2n-3}}
         \frac{f(x)}{F(x)^2(z_1 - x)(z_2- x)}\right]_{x=0}\\
      & - \frac{f(z_1)}{(z_1-z_2)z_1^{2n-2} F(z_1)^2}
          - \frac{f(z_2)}{(z_2-z_1) z_2^{2n-2}F(z_2)^2}
         .
\end{split}
\end{gather*}

\item
\begin{gather*}
  \frac{1}{z_1-x_i}\frac{1}{z_2-x_i}\left(
\frac{1}{z_1-x_i}+\frac{1}{z_2-x_i}\right)=
\frac{1}{z_2-z_1}\left(\frac{1}{(z_2-x_i)^2}
   -\frac{1}{(z_1-x_i)^2}\right)
         .
\end{gather*}

\item For symmetric term $K(i,j)=K(j,i)$,
\begin{gather*}
I:=\sum_{i\neq j} K(i,j)\frac{1}{z_1-x_i}
\frac{1}{z_2 -x_j}
\left(\frac{1}{x_j-x_i}-\frac{1}{z_2 - x_i}\right).
\end{gather*}
is expressed by
\begin{gather*}
I=\frac{1}{2}\sum_{i\neq j} \frac{K(i,j)}{z_1-z_2}
\left(\frac{1}{(z_1-x_j)(z_1-x_i)}
-\frac{1}{(z_2 - x_i)(z_2-x_j)}\right).
\end{gather*}
\end{enumerate}

\end{lemma}

\begin{proof}:
(1) will be proved by the following residual computation
Let $\partial C_g^o$ be the boundary of a polygon representation
$C_g^o$ of $C_g$,
\begin{gather*}
\oint_{\partial C_g^o} \omega =0,\quad
\omega:= \frac{f(x)}{(x-z_1)(x-z_2)x^{2n-2}F(x)^2}d x
 .
\end{gather*}
The divisor of the integrand is
\begin{gather*}
\begin{split}
\left(\omega\right)& =
     3\sum_{i=1}^{2g+1} (b_i,0) -
        2\sum_{i=1,\epsilon=\pm}^k (x_i,\epsilon y_i)
-\sum_{a=1,\epsilon=\pm}^2(z_a,\epsilon y(z_a)) \\
&\quad\quad\quad
-(2n-2)\sum_{\epsilon=\pm}(0,\epsilon y(0)) - (4g-4m+1)\infty .
\end{split}
\end{gather*}

\begin{enumerate}[{(}i{)}]
\item Using the fact that
the local parameter $t$ at $\infty$ is $x=1/t^2$, we have
\begin{gather*}
\res_{\infty}\omega
      = -2\frac{1}{(4g-4m)!}\left[
     \frac{\partial^{4g-4m}}{\partial t^{4g-4m}}
         \frac{t^{4g- 4k+1}f(1/t^2)}
        {F(1/t^2)^2(1 - z_1 t^2)(1- z_2 t^2)}\right]_{t=0}.
\end{gather*}

\item Since the local parameter at $x=0$ is $x$ itself, we have
\begin{gather*}
\res_{(0,\pm y(0))}\omega
      =  \frac{1}{(2n-1)!}
      \left[\frac{\partial^{2n-1}}{\partial x^{2n-1}}
         \frac{f(x)}{F(x)^2(z_1 - x)(z_2- x)}\right]_{x=0}.
\end{gather*}

\item Noting that the local parameter $t$ at $(x_k,\pm y_k)$ is
$t=x-x_k$, we have
\begin{gather*}
\res_{(x_k, \pm y_k)}\omega
      =  \frac{1}{F'(x_k)}
             \left[\frac{\partial}{\partial x}\left(
       \frac{f(x)}{(x - z_1)(x-z_2)x^{2n-2} F'(x) }\right) \right]_{x = x_k}
        .
\end{gather*}

\item Using the fact that the local parameter $t$ at $z_a$ is
$t=x-z_a$, we have
\begin{gather*}
\res_{(z_a,\pm y(z_a))}\omega
      =\frac{f(z_a)}{(z_a-z_b)z_a^{2n-2}F(z_a)^2}.
\end{gather*}
\end{enumerate}
where $z_b=z_2$ for $a=1$ and $z_1$ for $a=2$.

By arranging them, we obtain (1).

On the other hand, (2) can be proved by using a  trick:
\begin{gather*}
         \frac{1}{(z_1-x)(z_2-x)}
  =\frac{1}{z_1-z_2}
\left(\frac{1}{z_2-x}-\frac{1}{(z_1-x)}\right).
\end{gather*}
(3) can be evaluated by
$$
	I = 2\times\mbox{right hand side} - I.
$$
\end{proof}

\begin{proposition}
\begin{enumerate}

\item $m=g$ case:
$$
	H_\infty = 1.
$$

\item  $n=1$ case:
$$
	H_0 =0.
$$

\item  $n=2$ case:
$$
	H_0 =-\frac{2}{(x_1\cdots x_k)^2 z_1 z_2}\left[\lambda_1
        + \lambda_0\left(\sum_{i=1}^k \frac{2}{x_i}
          +\frac{1}{z_1} + \frac{1}{z_2}\right)\right].
$$
\end{enumerate}
\end{proposition}

\end{document}